\documentclass[11pt]{amsart}
\usepackage{latexsym}
\usepackage{amsmath}
\usepackage{amssymb}

\newcommand{\eproof}{\mbox{\ }\hfill $\Box$ \par \vskip 10pt}

\newtheorem{Theorem}{Theorem}[section]
\newtheorem{lemma}[Theorem]{Lemma}
\newtheorem{prop}[Theorem]{Proposition}

\numberwithin{equation}{section}

\def\cal{\mathcal}

\begin{document}

\title[Improved resolvent bounds]{Improved resolvent bounds for radial potentials. II}

\author[G. Vodev]{Georgi Vodev}

\address {Universit\'e de Nantes, Laboratoire de Math\'ematiques Jean Leray, 2 rue de la Houssini\`ere, BP 92208, 44322 Nantes Cedex 03, France}
\email{Georgi.Vodev@univ-nantes.fr}

\date{}

\begin{abstract} We prove semiclassical resolvent estimates for the
Schr\"odinger operator in $\mathbb{R}^d$, $d\ge 3$, with real-valued radial potentials $V\in L^\infty(\mathbb{R}^d)$.
We show that if $V(x)={\cal O}\left(\langle x\rangle^{-\delta}\right)$ with $\delta>4$, then the resolvent bound
is of the form $\exp\left(Ch^{-\frac{\delta}{\delta-1}}\left(\log(h^{-1})\right)^{\frac{1}{\delta-1}}\right)$ with some constant $C>0$. 
If $V(x)={\cal O}\left(e^{-\widetilde C\langle x\rangle^{\alpha}}\right)$ with $\widetilde C,\alpha>0$, we get better resolvent bounds
of the form $\exp\left(Ch^{-1}\left(\log(h^{-1})\right)^{\frac{1}{\alpha}}\right)$. 
\quad

Key words: Schr\"odinger operator, resolvent bounds, radial potentials.
\end{abstract} 

\maketitle

\setcounter{section}{0}
\section{Introduction and statement of results}

Our goal in this note is to improve some of the resolvent bounds proved in \cite{kn:V4} for radial real-valued potentials. We also give 
a different, shorter proof of the sharp resolvent bounds proved recently in \cite{kn:DGS} for radial compactly supported real-valued potentials.
Consider the Schr\"odinger operator
$$P(h)=-h^2\Delta+V(x)$$
where $0<h\ll 1$ is a semiclassical parameter, $\Delta$ is the negative Laplacian in 
$\mathbb{R}^d$, $d\ge 3$, and $V\in L^\infty(\mathbb{R}^d)$ is a real-valued short-range potential satisfying the condition
\begin{equation}\label{eq:1.1}
|V(x)|\le C(|x|+1)^{-\delta}
\end{equation}
where $C>0$ and $\delta>1$ are some constants. We are interested in bounding  
the quantity
$$g_s^\pm(h,\varepsilon):=\log\left\|(|x|+1)^{-s}(P(h)-E\pm i\varepsilon)^{-1}(|x|+1)^{-s}
\right\|_{L^2(\mathbb{R}^d)\to L^2(\mathbb{R}^d)}$$
from above by an explicit function of $h$, independent of $\varepsilon$. 
Here $0<\varepsilon<1$, $s>1/2$ is independent of $h$ and $E>0$ is a fixed energy level independent of $h$.
When $\delta>2$ it has been proved in \cite{kn:GS1} that
\begin{equation}\label{eq:1.2}
g_s^\pm(h,\varepsilon)\le Ch^{-4/3}\log(h^{-1}).
\end{equation}
The bound (\ref{eq:1.2}) was first proved in \cite{kn:KV} and \cite{kn:S} for compactly supported potentials,
and in \cite{kn:V1} when $\delta>3$. It was also shown in \cite{kn:V2} that (\ref{eq:1.2}) still holds for
more general asymptotically Euclidean manifolds. On the other hand, it is shown in \cite{kn:V4} 
that the logarithmic term in the righ-hand side of
(\ref{eq:1.2}) can be removed for real-valued potentials $V$ depending only on the radial variable $r=|x|$, 
provided $V$ satisfies (\ref{eq:1.1}) with $\delta>2$. Furhtemore, for compactly supported radial potentials
a much better bound has been recently proved in \cite{kn:DGS}, namely the following one
\begin{equation}\label{eq:1.3}
g_s^\pm(h,\varepsilon)\le Ch^{-1}.
\end{equation}
The bound (\ref{eq:1.3}) was previously proved in \cite{kn:D}, \cite{kn:GS2}, \cite{kn:V3} for slowly decaying Lipschitz potentials $V$ 
with respect to the radial variable $r$. Note also that when $d=1$ the bound (\ref{eq:1.3}) is proved in 
 \cite{kn:DS} for $V\in L^1(\mathbb{R})$. 
We show in the present paper that better resolvent bounds than those obtained 
in \cite{kn:V4} can be proved for 
non-compactly supported radial potentials, too. 
Our main result is the following

\begin{Theorem} Let $d\ge 3$ and suppose that the potential $V$ depends only on the radial variable. If $V$ 
satisfies (\ref{eq:1.1}) with $\delta>4$,  
then there exist constants $C>0$ and $h_0>0$ independent of $h$ and $\varepsilon$ but depending on $s$, $E$, such 
that the bound
\begin{equation}\label{eq:1.4}
g_s^\pm(h,\varepsilon)\le Ch^{-\frac{\delta}{\delta-1}}\left(\log(h^{-1})\right)^{\frac{1}{\delta-1}}
\end{equation}
holds for all $0<h\le h_0$. If the potential satisfies the condition
\begin{equation}\label{eq:1.5}
|V(x)|\le C_1e^{-C_2|x|^{\alpha}}
\end{equation}
with some constants $C_1,C_2,\alpha>0$, then we have the better bound
\begin{equation}\label{eq:1.6}
g_s^\pm(h,\varepsilon)\le Ch^{-1}\left(\log(h^{-1})\right)^{\frac{1}{\alpha}}.
\end{equation}
\end{Theorem}

Note that when $V$ is compactly supported the proof of Theorem 1.1 leads to the bound (\ref{eq:1.3}) already proved in \cite{kn:DGS} 
in a different way. 

The fact that the potential is radial plays an important role in the proof of the above theorem. It allows us to reduce
the $d$ - dimensional resolvent bound to one-dimensional ones. In other words, we have to 
bound the resolvent of an infinite family of one-dimensional Schr\"odinger operators depending on an additional parameter
denoted by $\nu\ge 0$ below, which can be expressed in terms of the eigenvalues of the
Laplace-Beltrami operator on the $d-1$ - dimensional unit sphere (see Section 2). To do so,
we make use of some bounds already proved in \cite{kn:V4} (see Proposition 2.2) 
and we show that Proposition 3.2 of \cite{kn:V4} can be improved
significantly for potentials decaying at infinity sufficiently fast (see Proposition 2.3). 
It is not clear if the bounds in Theorem 1.1 still hold for non-radial $L^\infty$ potentials
since neither a proof nor counterexamples are available. To author's best knoweldge, the best resolvent bound
for such potentials is (\ref{eq:1.2}), which seems hard to improve without extra conditions even if the potential is
supposed to be compactly supported.

\section{Preliminaries}

We will use the fact that the potential is radial to reduce the resolvent bound to infinitely many one-dimensional 
resolvent bounds (see also Section 2 of \cite{kn:V4}) .
To this end we will write the operator $P(h)$ in 
polar coordinates $(r,w)\in\mathbb{R}^+\times\mathbb{S}^{d-1}$, $r=|x|$, $w=x/|x|$ and we will use that $L^2(\mathbb{R}^d)=L^2(\mathbb{R}^+\times\mathbb{S}^{d-1}, r^{d-1}drdw)$. We have the identity
\begin{equation}\label{eq:2.1}
 r^{(d-1)/2}\Delta  r^{-(d-1)/2}=\partial_r^2+\frac{\widetilde\Delta_w}{r^2}
\end{equation}
where $\widetilde\Delta_w=\Delta_w-\frac{1}{4}(d-1)(d-3)$ and $\Delta_w$ denotes the negative Laplace-Beltrami operator
on $\mathbb{S}^{d-1}$. Using (\ref{eq:2.1}) we can write the operator
$${\cal P}^\pm(h)=r^{(d-1)/2}(P(h)-E\pm i\varepsilon)r^{-(d-1)/2}$$
  in the coordinates $(r,w)$ as follows
$${\cal P}^\pm(h)={\cal D}_r^2+\frac{\Lambda_w}{r^2}+V(r)-E\pm i\varepsilon $$
where we have put ${\cal D}_r=-ih\partial_r$ and $\Lambda_w=-h^2\widetilde\Delta_w$. 
Let $\lambda_j\ge 0$ be the eigenvalues of $-\Delta_w$ repeated with the multiplicities and let $e_j\in L^2(\mathbb{S}^{d-1})$
be the corresponding eigenfunctions. Set
$$\nu=h\sqrt{\lambda_j+\frac{1}{4}(d-1)(d-3)}$$
and
$$Q^\pm_\nu(h)={\cal D}_r^2+\frac{\nu^2}{r^2}+V(r)-E\pm i\varepsilon.$$
Let $v\in L^2(\mathbb{R}^+\times\mathbb{S}^{d-1},drdw)$ and set 
$$v_j(r)=\langle v(r,\cdot),e_j\rangle_{L^2(\mathbb{S}^{d-1})}.$$ 
We can write 
$$v=\sum_j v_je_j,$$
$${\cal P}^\pm(h)v=\sum_jQ^\pm_\nu(h)v_j e_j,$$
so we have 
$$\|v\|^2_{L^2(\mathbb{R}^+\times\mathbb{S}^{d-1})}=\sum_j\|v_j\|^2_{L^2(\mathbb{R}^+)},$$
$$\|(r+1)^{-s}v\|^2_{L^2(\mathbb{R}^+\times\mathbb{S}^{d-1})}=\sum_j\|(r+1)^{-s}v_j\|^2_{L^2(\mathbb{R}^+)},$$
$$\|(r+1)^{s}{\cal P}^\pm(h)v\|^2_{L^2(\mathbb{R}^+\times\mathbb{S}^{d-1})}=\sum_j\|(r+1)^{s}Q_\nu^\pm(h)v_j\|^2_{L^2(\mathbb{R}^+)}.$$
The following lemma is proved in Section 2 of \cite{kn:V4} using the above identities.

\begin{lemma} Let $s>1/2$ and suppose that for all $\nu$ the estimates
\begin{equation}\label{eq:2.2}
\|(r+1)^{-s}u\|^2_{L^2(\mathbb{R}^+)}
\le M_\nu\|(r+1)^{s}Q_\nu^\pm(h)u\|^2_{L^2(\mathbb{R}^+)}$$
$$+M_\nu\varepsilon\|u\|^2_{L^2(\mathbb{R}^+)}
+M_\nu\varepsilon\|{\cal D}_ru\|^2_{L^2(\mathbb{R}^+)}
\end{equation}
hold for every $u\in H^2(\mathbb{R}^+)$
such that $u(0)=0$ and $(r+1)^{s}Q_\nu^\pm(h)u\in L^2(\mathbb{R}^+)$, with $M_\nu>0$ independent of $\varepsilon$ and $u$.  
Suppose also that 
$$M:=\left(2+E+\|V\|_{L^\infty}\right)\sup_{\nu^2\in{\rm spec}\,\Lambda_w}M_\nu<\infty.$$
Then we have the bound
\begin{equation}\label{eq:2.3}
g_s^\pm(h,\varepsilon)\le \log(M+1).
\end{equation}
\end{lemma}

Thus we reduce our problem to proving estimates like (\ref{eq:2.2}) with as good bounds $M_\nu$ as possible. 
We will make use of the following proposition proved in Section 3 of \cite{kn:V4} (see Proposition 3.1 of \cite{kn:V4}).

\begin{prop} The estimate (\ref{eq:2.2}) holds for all 
$\nu$ with $M_\nu=e^{C(\nu+1)/h}$, where $C>0$ is a constant independent of $\nu$ and $h$.
\end{prop} 

Therefore, we only need to bound $M_\nu$ for large $\nu$. Set $\tau=1$ if $V$ is compactly supported,
$\tau=(\epsilon h)^{-\frac{1}{\delta-1}}$ if $V$ satisfies (\ref{eq:1.1}) and $\tau=\epsilon^{-1/\alpha}$ 
if $V$ satisfies (\ref{eq:1.5}),
where $\epsilon=(\log(h^{-1}))^{-1}\ll 1$.  
In what follows in this paper we will prove the following

\begin{prop} The estimate (\ref{eq:2.2}) holds for all 
$\nu\ge c\tau$ with $M_\nu=C(\epsilon h)^{-2}$, where $C,c>0$ are constants independent of $\nu$ and $h$.
\end{prop} 

Clearly, the bounds (\ref{eq:1.4}) and (\ref{eq:1.6}) follow from (\ref{eq:2.3}) and Propositions 2.2 and 2.3.

\section{A priori estimates}

Let $\phi_0\in C^\infty(\mathbb{R})$ be a
real-valued function such that $0\le\phi_0\le 1$, $\phi'_0\ge 0$, $\phi_0(\sigma)=0$ for $\sigma\le 1$, 
$\phi_0(\sigma)=1$ for $\sigma\ge 2$, and set $\phi(r)=\phi_0(r/\lambda)$, where $\lambda\gg 1$. We also set 
$$Q^\pm_{\nu,0}(h)={\cal D}_r^2+\frac{\nu^2}{r^2}-E\pm i\varepsilon.$$
In this section we will prove the following

\begin{prop} For every $u\in H^2(\mathbb{R}^+)$
such that $(r+1)^{1+\epsilon}Q_{\nu,0}^\pm(h)u\in L^2(\mathbb{R}^+)$, we have the estimate
\begin{equation}\label{eq:3.1}
\int_0^\infty(r+1)^{-1-\epsilon}\left(|\phi u(r)|^2+|{\cal D}_r(\phi u)(r)|^2\right)dr$$ 
$$\le C\lambda^{-1}\epsilon^{-1}h\int_{\lambda}^{2\lambda}\left(|u(r)|^2+|{\cal D}_ru(r)|^2\right)dr$$ 
$$+C(\epsilon h)^{-2}\int_0^\infty(r+1)^{1+\epsilon}|\phi Q_{\nu,0}^\pm(h)u(r)|^2dr$$
$$+C(\epsilon h)^{-1}\varepsilon\int_0^\infty\left(|u(r)|^2+|{\cal D}_ru(r)|^2\right)dr
\end{equation}
with a constant $C>0$ independent of $\epsilon$, $\varepsilon$, $\nu$, $\lambda$ and $h$.
\end{prop}

{\it Proof.} It is easy to see that the first derivative of the function 
$$F(r)=(E-\nu^2r^{-2})|\phi u(r)|^2+|{\cal D}_r(\phi u)(r)|^2.$$
 is given by 
$$F'(r)=2\nu^2r^{-3}|\phi u|^2-\Phi(r)$$
where
$$\Phi(r)=2h^{-1}{\rm Im}\,Q_{\nu,0}^\pm(h)\phi u\overline{{\cal D}_r(\phi u)}\mp 
2\varepsilon h^{-1}{\rm Re}\,\phi u\overline{{\cal D}_r(\phi u)}$$
$$=2h^{-1}{\rm Im}\,\phi Q_{\nu,0}^\pm(h)u\overline{{\cal D}_r(\phi u)}\mp 
2\varepsilon h^{-1}{\rm Re}\,\phi u\overline{{\cal D}_r(\phi u)}+\Psi(r)$$
$$\le \gamma(r+1)^{-1-\epsilon}\left|{\cal D}_r(\phi u)\right|^2+\gamma^{-1}h^{-2}(r+1)^{1+\epsilon}\left|\phi Q_{\nu,0}^\pm(h)u\right|^2$$
$$+\varepsilon h^{-1}\left(|u|^2+|{\cal D}_ru|^2\right)+\Psi(r),$$
$\gamma>0$ being arbitrary, where 
$$\Psi(r)=2h^{-1}{\rm Im}\,[{\cal D}_r^2,\phi]u\overline{{\cal D}_r(\phi u)}$$
$$=-2{\rm Im}\left(2i\phi'{\cal D}_ru+h\phi''u\right)\left(\phi\overline{{\cal D}_ru}+ih\phi'\overline{u}\right)$$
$$=-4\phi\phi'|{\cal D}_ru|^2+2h(2\phi'^2+\phi\phi''){\rm Im}\,\overline{u}{\cal D}_ru-2h^2\phi'\phi''|u|^2$$
$$\le Ch\lambda^{-2}\left(\phi'_0(r/\lambda)+|\phi''_0(r/\lambda)|\right)\left(|u(r)|^2+|{\cal D}_ru(r)|^2\right)$$
with some constant $C>0$ independent of $h$ and $\lambda$. Hence, given any $t>0$, we get
$$F(t)=-\int_t^\infty F'(r)dr\le \int_0^\infty\Phi(r)dr$$
$$\le \gamma\int_0^\infty(r+1)^{-1-\epsilon}\left|{\cal D}_r(\phi u)\right|^2dr$$
$$+\gamma^{-1}h^{-2}\int_0^\infty(r+1)^{1+\epsilon}\left|\phi Q_{\nu,0}^\pm(h)u\right|^2dr$$
$$+h^{-1}\varepsilon\int_0^\infty\left(|u|^2+|{\cal D}_ru|^2\right)dr$$ $$
+Ch\lambda^{-2}\int_{\lambda}^{2\lambda}\left(|u|^2+|{\cal D}_ru|^2\right)dr.$$
Multiplying this inequality by $(t+1)^{-1-\epsilon}$ and integrating with respect to $t$ lead to the estimate
 $$\int_0^\infty (t+1)^{-1-\epsilon}F(t)dt\le \gamma\int_0^\infty(r+1)^{-1-\epsilon}\left|{\cal D}_r(\phi u)\right|^2dr$$
$$+\gamma^{-1}(\epsilon h)^{-2}\int_0^\infty(r+1)^{1+\epsilon}\left|\phi Q_{\nu,0}^\pm(h)u\right|^2dr$$
$$+(\epsilon h)^{-1}\varepsilon\int_0^\infty\left(|u|^2+|{\cal D}_ru|^2\right)dr$$ $$
+Ch\lambda^{-2}\epsilon^{-1}\int_{\lambda}^{2\lambda}\left(|u|^2+|{\cal D}_ru|^2\right)dr$$
for any $\gamma>0$, where we have used that 
$$\int_0^\infty (t+1)^{-1-\epsilon}dt=\epsilon^{-1}.$$
Taking $\gamma$ small enough, indpendent of $\epsilon$, $h$ and $\lambda$, we can absorb the first term in the right-hand side
of the above inequality. Thus we get
$$ \int_0^\infty(r+1)^{-1-\epsilon}\left(|\phi u(r)|^2+|{\cal D}_r(\phi u)(r)|^2\right)dr$$ 
$$\lesssim \int_0^\infty\nu^2r^{-3}\left|\phi u\right|^2dr$$
$$+(\epsilon h)^{-2}\int_0^\infty(r+1)^{1+\epsilon}\left|\phi Q_{\nu,0}^\pm(h)u\right|^2dr$$
$$+(\epsilon h)^{-1}\varepsilon\int_0^\infty\left(|u|^2+|{\cal D}_ru|^2\right)dr$$ $$
+h\lambda^{-2}\epsilon^{-1}\int_{\lambda}^{2\lambda}\left(|u|^2+|{\cal D}_ru|^2\right)dr.$$
On the other hand, we have
$$\int_0^\infty\nu^2r^{-3}\left|\phi u\right|^2dr=\int_0^\infty F'(r)dr+\int_0^\infty\Phi(r)dr=\int_0^\infty\Phi(r)dr$$
$$\le \gamma\int_0^\infty(r+1)^{-1-\epsilon}\left|{\cal D}_r(\phi u)\right|^2dr$$
$$+\gamma^{-1}h^{-2}\int_0^\infty(r+1)^{1+\epsilon}\left|\phi Q_{\nu,0}^\pm(h)u\right|^2dr$$
$$+h^{-1}\varepsilon\int_0^\infty\left(|u|^2+|{\cal D}_ru|^2\right)dr$$ $$
+Ch\lambda^{-2}\int_{\lambda}^{2\lambda}\left(|u|^2+|{\cal D}_ru|^2\right)dr$$
for any $\gamma>0$. Combining the above inequalities and taking 
$\gamma$ small enough, indpendent of $\epsilon$, $h$ and $\lambda$, in order to absorb the corresponding term, we get (\ref{eq:3.1}).
\eproof

\section{Proof of Proposition 2.3}

We will first derive from Proposition 3.1 the following 

\begin{prop} Let $u\in H^2(\mathbb{R}^+)$ be 
such that $(r+1)^{1+\epsilon}Q_\nu^\pm(h)u\in L^2(\mathbb{R}^+)$. There exists a constant $\lambda_0>0$ such that
if $\lambda\ge\lambda_0\tau$, then we have the estimate
\begin{equation}\label{eq:4.1}
\int_{2\lambda}^\infty(r+1)^{-1-\epsilon}\left(|u(r)|^2+|{\cal D}_ru(r)|^2\right)dr$$ 
$$\le C\lambda^{-1}\epsilon^{-1}h\int_{\lambda}^{2\lambda}\left(|u(r)|^2+|{\cal D}_ru(r)|^2\right)dr$$ 
$$+C(\epsilon h)^{-2}\int_0^\infty(r+1)^{1+\epsilon}|Q_\nu^\pm(h)u(r)|^2dr$$
$$+C(\epsilon h)^{-1}\varepsilon\int_0^\infty\left(|u(r)|^2+|{\cal D}_ru(r)|^2\right)dr
\end{equation}
with a constant $C>0$ independent of $\epsilon$, $\varepsilon$, $\nu$, $\lambda$ and $h$.
\end{prop}

{\it Proof.} We apply the estimate (\ref{eq:3.1}) and observe that
\begin{equation}\label{eq:4.2}
\int_0^\infty(r+1)^{1+\epsilon}|\phi Q_{\nu,0}^\pm(h)u(r)|^2dr$$
$$\le \int_0^\infty(r+1)^{1+\epsilon}|\phi Q_{\nu}^\pm(h)u(r)|^2dr$$ $$+
\int_0^\infty(r+1)^{1+\epsilon}|V(r)(\phi u)(r)|^2dr$$
$$\le \int_0^\infty(r+1)^{1+\epsilon}|Q_{\nu}^\pm(h)u(r)|^2dr$$ $$+
\rho(\lambda)\int_0^\infty(r+1)^{-1-\epsilon}|(\phi u)(r)|^2dr
\end{equation}
where
$$\rho(\lambda)=\sup_{r\ge\lambda}\,(r+1)^{2+2\epsilon}|V(r)|^2.$$
When $V$ is compactly supported we have $\rho(\lambda)=0$, provided $\lambda=\lambda_0$ is big enough, independent of $h$.
When $V$ satisfies (\ref{eq:1.1}) with $\delta>4$ we have
$$\rho(\lambda)\lesssim \lambda^{-2\delta+2+2\epsilon}\lesssim \lambda_0^{-2\delta+2}(\epsilon h)^2.$$
When $V$ satisfies (\ref{eq:1.5}) we have
$$\rho(\lambda)\lesssim \lambda^{2+2\epsilon}e^{-C_2\lambda^\alpha}
\lesssim \lambda_0^2\epsilon^{-2/\alpha}e^{C_2\lambda_0^\alpha\log h}\lesssim \epsilon^{-2/\alpha}h^{C_2\lambda_0^\alpha}
\lesssim h(\epsilon h)^2$$
provided $\lambda_0$ is big enough, independent of $h$ and $\epsilon$. Thus, taking $h$ small enough or 
$\lambda_0$ big enough, we can absorb the last term in the right-hand
side of (\ref{eq:4.2}) and obtain (\ref{eq:4.1}).
\eproof

To prove Proposition 2.3 we will combine Proposition 4.1 with the following

\begin{prop} Let $u\in H^2(\mathbb{R}^+)$ be 
such that $u(0)=0$. Then there exists a constant $\kappa>0$ such that we have the estimate
\begin{equation}\label{eq:4.3}
\int_0^{3\kappa\nu}\left(|u(r)|^2+|{\cal D}_ru(r)|^2\right)dr$$ 
$$\le Ch^2\nu^{-2}\int_{3\kappa\nu}^{4\kappa\nu}\left(|u(r)|^2+|{\cal D}_ru(r)|^2\right)dr+4\int_0^{\infty}|Q_\nu^\pm(h)u(r)|^2dr
\end{equation}
with a constant $C>0$ independent of $\varepsilon$, $\nu$ and $h$.
\end{prop}

{\it Proof.}
Let $\psi_0\in C^\infty(\mathbb{R})$ be a
real-valued function such that $0\le\psi_0\le 1$, $\psi_0(\sigma)=1$ for $\sigma\le 3$, 
$\psi_0(\sigma)=0$ for $\sigma\ge 4$, and set $\psi(r)=\psi_0(r/\kappa\nu)$, where $\kappa^{-1}=4\sqrt{1+E+\|V\|_{L^\infty}}$.
The estimate (\ref{eq:4.3}) is a consequence of the following

\begin{lemma} We have the estimate
\begin{equation}\label{eq:4.4}
\int_0^\infty\left(|\psi u(r)|^2+|{\cal D}_r(\psi u)(r)|^2\right)dr\le 4\int_0^{\infty}|Q_\nu^\pm(h)(\psi u)(r)|^2dr.
\end{equation}
\end{lemma}

{\it Proof.} The choice of $\kappa$ guarantees the inequality
$$(\nu^2r^{-2}+V(r)-E)|\psi u(r)|^2\ge |\psi u(r)|^2$$
for all $r$. 
Therefore, integrating by parts we obtain
$${\rm Re}\int_0^\infty Q_\nu^\pm(h)\psi u(r)\overline{\psi u(r)}dr$$
$$=\int_0^\infty|{\cal D}_r(\psi u)(r)|^2
+\int_0^\infty(\nu^2r^{-2}+V(r)-E)|\psi u(r)|^2dr$$
$$\ge \int_0^\infty|{\cal D}_r(\psi u)(r)|^2dr+\int_0^\infty |\psi u(r)|^2dr.$$
Hence
$$\int_0^\infty|{\cal D}_r(\psi u)(r)|^2dr+\int_0^\infty |\psi u(r)|^2dr$$
$$\le 2\int_0^{\infty}|Q_\nu^\pm(h)(\psi u)(r)|^2dr+
\frac{1}{2}\int_0^{\infty}|\psi u(r)|^2dr,$$
which clearly implies (\ref{eq:4.4}).
\eproof

Since 
$$\int_0^{\infty}|[Q_\nu^\pm(h),\psi]u(r)|^2dr=\int_0^{\infty}|[{\cal D}_r^2,\psi]u(r)|^2dr$$
$$\le Ch^2\nu^{-2}\int_{3\kappa\nu}^{4\kappa\nu}\left(|u(r)|^2+|{\cal D}_ru(r)|^2\right)dr,$$
the estimate (\ref{eq:4.3}) follows from (\ref{eq:4.4}).
\eproof

Let $u\in H^2(\mathbb{R}^+)$ be 
such that $u(0)=0$ and $(r+1)^{1+\epsilon}Q_\nu^\pm(h)u\in L^2(\mathbb{R}^+)$. 
We apply Proposition 4.1 with $\lambda=\kappa\nu$ and suppose that $\nu\ge \lambda_0\tau/\kappa$. By 
(\ref{eq:4.1}) and (\ref{eq:4.3}) we have
$$(4\kappa\nu+1)^{-1-\epsilon}\int_{3\kappa\nu}^{4\kappa\nu}\left(|u(r)|^2+|{\cal D}_ru(r)|^2\right)dr$$ 
$$\le \int_{2\kappa\nu}^\infty(r+1)^{-1-\epsilon}\left(|u(r)|^2+|{\cal D}_ru(r)|^2\right)dr$$ 
$$\le C\nu^{-1}\epsilon^{-1}h\int_{\kappa\nu}^{2\kappa\nu}\left(|u(r)|^2+|{\cal D}_ru(r)|^2\right)dr$$ 
$$+C(\epsilon h)^{-2}\int_0^\infty(r+1)^{1+\epsilon}|Q_\nu^\pm(h)u(r)|^2dr$$
$$+C(\epsilon h)^{-1}\varepsilon\int_0^\infty\left(|u(r)|^2+|{\cal D}_ru(r)|^2\right)dr$$
$$\le C\nu^{-3}\epsilon^{-1}h^3\int_{3\kappa\nu}^{4\kappa\nu}\left(|u(r)|^2+|{\cal D}_ru(r)|^2\right)dr$$ 
$$+C(\epsilon h)^{-2}\int_0^\infty(r+1)^{1+\epsilon}|Q_\nu^\pm(h)u(r)|^2dr$$
$$+C(\epsilon h)^{-1}\varepsilon\int_0^\infty\left(|u(r)|^2+|{\cal D}_ru(r)|^2\right)dr.$$
It is clear that taking $h$ small enough we can absorb the first term in the right-hand side of the above inequality. Thus we obtain
$$(4\kappa\nu+1)^{-1-\epsilon}\int_{3\kappa\nu}^{4\kappa\nu}\left(|u(r)|^2+|{\cal D}_ru(r)|^2\right)dr$$ 
$$\le C(\epsilon h)^{-2}\int_0^\infty(r+1)^{1+\epsilon}|Q_\nu^\pm(h)u(r)|^2dr$$
$$+C(\epsilon h)^{-1}\varepsilon\int_0^\infty\left(|u(r)|^2+|{\cal D}_ru(r)|^2\right)dr,$$
which together with (\ref{eq:4.3}) yield
\begin{equation}\label{eq:4.5}
\int_0^{3\kappa\nu}\left(|u(r)|^2+|{\cal D}_ru(r)|^2\right)dr$$ 
$$\le C\epsilon^{-2}\int_0^\infty(r+1)^{1+\epsilon}|Q_\nu^\pm(h)u(r)|^2dr$$
$$+C\varepsilon\int_0^\infty\left(|u(r)|^2+|{\cal D}_ru(r)|^2\right)dr
\end{equation}
with a new constant $C>0$. Combining (\ref{eq:4.1}) with (\ref{eq:4.5}) leads to the estimate
\begin{equation}\label{eq:4.6}
\int_{2\kappa\nu}^\infty(r+1)^{-1-\epsilon}\left(|u(r)|^2+|{\cal D}_ru(r)|^2\right)dr$$ 
$$\le C(\epsilon h)^{-2}\int_0^\infty(r+1)^{1+\epsilon}|Q_\nu^\pm(h)u(r)|^2dr$$
$$+C(\epsilon h)^{-1}\varepsilon\int_0^\infty\left(|u(r)|^2+|{\cal D}_ru(r)|^2\right)dr.
\end{equation}
By (\ref{eq:4.5}) and (\ref{eq:4.6}) we conclude
\begin{equation}\label{eq:4.7}
\int_0^\infty(r+1)^{-1-\epsilon}\left(|u(r)|^2+|{\cal D}_ru(r)|^2\right)dr$$ 
$$\le C(\epsilon h)^{-2}\int_0^\infty(r+1)^{1+\epsilon}|Q_\nu^\pm(h)u(r)|^2dr$$
$$+C(\epsilon h)^{-1}\varepsilon\int_0^\infty\left(|u(r)|^2+|{\cal D}_ru(r)|^2\right)dr.
\end{equation}
Taking $h$ small enough we can arrange that $\epsilon<2s-1$.  Therefore, Proposition 2.3 follows from (\ref{eq:4.7}).
\eproof

\end{document}